\documentclass[french,12pt,ams]{article}
\usepackage{amssymb,amsfonts,amsmath}
\usepackage[french]{babel} 
\usepackage[latin1]{inputenc}
\usepackage{graphicx}
\DeclareGraphicsExtensions{.pdf,.jpg,.png, .bmp, .eps}

\newtheorem{lem}{Lemma}[section]

\newtheorem{theo}{Theorem}

\newcommand{\N}{{\bf N}}
\newcommand{\R}{{\bf R}}
\newcommand{\Z}{{\bf Z}}

\newenvironment{nproof}[1]{\trivlist\item[\hskip \labelsep{\bf Proof{#1}.}]}
{\begin{flushright} $\square$\end{flushright}\endtrivlist}

\newenvironment{block}[1]{\trivlist\item[\hskip \labelsep{{#1}.}]}{\endtrivlist}

\newtheorem{com}{Comment}

\newcommand{\Ci}{\bf{C}}

{\end{enumerate}\end{sc}}

\title{
On the Hausdorff dimension of Julia sets of some real polynomials}
\author{Genadi Levin, Michel Zinsmeister}

\begin{document}
\maketitle
\thanks{Supported in part by IMPAN-BC European Community
Centre of Excellence and by the Marie Curie European network CODY}

\abstract{ We show that the supremum for $c$ real of the Hausdorff dimension 
of the Julia set of the polynomial $z\mapsto z^d+c$ ($d$ is an even natural number) 
is greater than $2d/(d+1)$. }
\section{Introduction and statement of the result}\label{s1}
The Julia set of a non-linear polynomial $P: \Ci\to \Ci$ 
is the set of points  having no 
neighborhood on which the family of iterates $(P^n)$ is normal.
It is a compact non-empty set which is, except for very special polynomials 
$P$, a fractal set.
For $c\in \Ci$ we denote by $f_c$ the polynomial
$$f_c(z)=z^d+c,$$
where $d\ge 2$ is an even integer number, which we fix. 
Denote by $J_c$ the Julia set of $f_c$.

In the quadratic case $d=2$, the values of $c$ for which the Hausdorff
dimension (HD) of $J_c$ is big has attracted a lot of attention.
It has been known since the pioneering work of Douady-Hubbard ~\cite{DH}
that the Hausdorff dimension of the Julia set
is less than $2$ for every hyperbolic polynomial. Thus
$HD(J_c)<2$ outside $\partial M$ and the (hypothetic) non-hyperbolic
components of (the interior) $M$. We recall that $M$ stands for the
Mandelbrot set, that is the compact subset of $c\in \Ci$ such that
$J_c$ is connected. Shishikura ~\cite{Shi} was the first to find quadratic
Julia sets with Hausdorff dimension $2$. He indeed proved that this property
holds on a dense $\cal{G}_\delta$ subset of $\partial M$ or even
on a dense $\cal{G}_\delta$ subset on the boundary of every
hyperbolic component of $M$. More recently, Buff and Ch\'eritat ~\cite{BC}
have found quadratic Julia sets with positive Lebesgue measure (see also http://annals.math.princeton.edu/articles/3682).

Both Shishikura's and Buff-Ch\'eritat's results are based on the phenomenon
of parabolic implosion which has been discovered and studied by 
Douady-Hubbard \cite{DH}. It should be pointed out that 
Buff-Ch\'eritat's result is very involved and that we will make no use of it.
There is no doubt that if they exist, values of $c\in \bf{R}$
such that the Julia set $J_c$ has positive measure must be as hard to find
as Buff-Cheritat's ones. The aim of present note is much more modest.
Its starting point is the second author re-readind of Shishikura's 
result ~\cite{Z}: it states that if one implodes a polynomial
with a parabolic cycle having $q$ petals, then the dimension
of its Julia set
automatically gets bigger than $2q/(q+1)$. Shishikura's result 
follows from this by a Baire argument. 

Very little is known about Hausdorff dimension of $J_c$ for real $c$.
In particular, it is not known if, for a given degree $d$,
$$\sup\{HD(J_c), c\in \bf {R}\}=2.$$
Possible candidates for high dimension are of course (just look at them!)
infinitely renormalizable polynomials but the analysis seems to be 
very delicate and at least no result concerning dimension $2$
has been proven so far (for results 
in the opposite direction, see ~\cite{AL} though).
It is for example unknown if the Julia set of
the Feigenbaum polynomial has dimension $2$
or not (see~\cite{LS1}-~\cite{LS2} for the Julia set of
the Feigenbaum universal map though). The only known result
about this set is a very general result of Zdunik \cite{Zd}: it has
dimension bigger than $1$.

If one tries to use the same ideas as in ~\cite{Shi} for real polynomials
$f_c$, one immediately faces the problem that if $f_c$, for $c$ real, 
has a parabolic cycle that may be imploded along the real axis
then the number of petals is $1$ and this does not imply more than Zdunik's 
general result. The only trick of this paper is to make use of a virtual 
doubling of petals when the critical point is mapped to a parabolic
point (by Lavaurs map). It was inspired by Douady et al paper~\cite{BDDS} and
implies the following theorem, which is the main result of this work.
\begin{theo}\label{t}
Let $f_c(z)=z^d+c$, $d$ even.
Let $N$ be the set of parameters $c\in \bf{R}$, such that
$f_c$ has a parabolic cycle of period at least $2$ and multiplier $1$.
Then there exists an open set $Y$ of $\bf{R}$ whose closure contains $N$ 
such that $J_c$ is connected and 
$$HD(J_c)>\frac{2d}{d+1},$$ 
for every $c\in Y$.
\end{theo}
\begin{com}
In fact, we prove a stronger statement: hyperbolic dimension~\cite{Shi}
of $J_c$ is bigger than $2d/(d+1)$.
\end{com}
\begin{com}
By ~\cite{KSS} (~\cite{GS},~\cite{LY} for $d=2$),
the set of real $c$ such that $f_c$ is hyperbolic is dense in $\bf{R}$.
In particular, 
hyperbolic parameters $c$ are dense in $Y$. 
\end{com}
{\bf Acknowledgment.} 
This work was done during the first author's one month's visit
at the university of Orl\'eans in 2008.
 \section{Proof of the theorem}
We fix an even integer $d\ge 2$ and consider the family $f_c(x)=x^d+c$,
for real $c$.
Then the Julia set $J_c$ is connected
if and only if $c\in [a, b]$, where $a<0$ is such that
$f_a^2(0)$ is a fixed point, and $b>0$ is such that $f_b$
has a fixed point with multiplier $1$.
It is sufficient to prove 
that, given $c_0\in (a,b)$ such that $f_{c_0}$ has a neutral cycle 
of period $k>1$ with multiplier $1$, that is, with
one petal, there is an open set $Z$ accumulating
at $c_0$ for which $HD(J_c)>2d/(d+1)$ for $c\in\Z$.

We begin with three pictures. The first two ones  illustrate 
how to choose the parameter (for $d=2$) and how the corresponding Julia set looks like. The third one has been kindly drawn for us by the referee and shows the case $d=4$.
\begin{figure}[!ht]
\begin{center}
\includegraphics[width=10cm,angle=270]{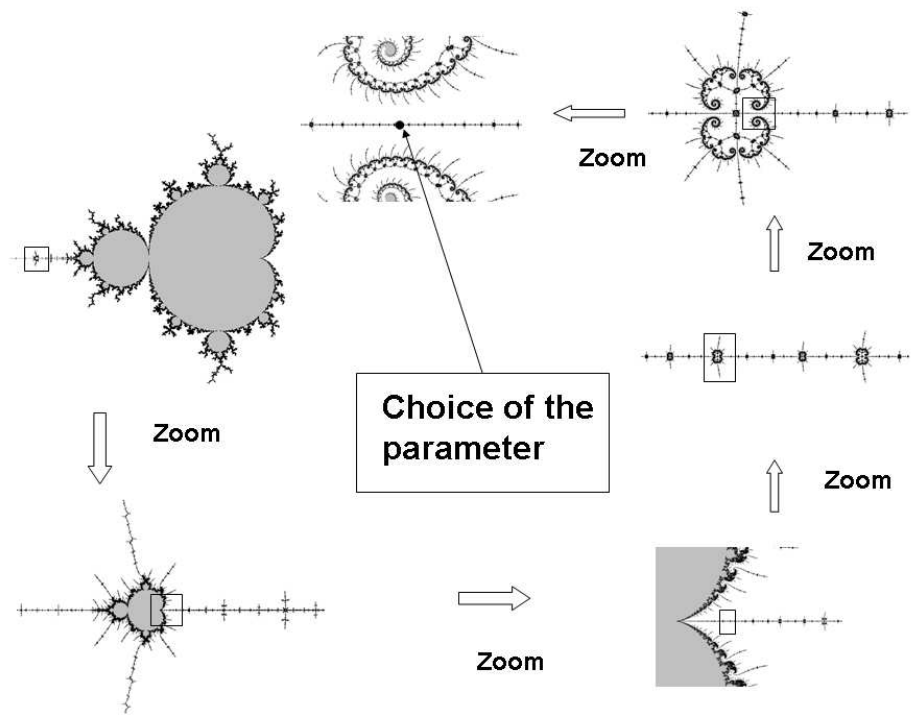}
\caption{Choice of the parameter}
\end{center}
\end{figure}
\begin{figure}[!ht]
\begin{center}
\includegraphics[width=10cm,angle=270]{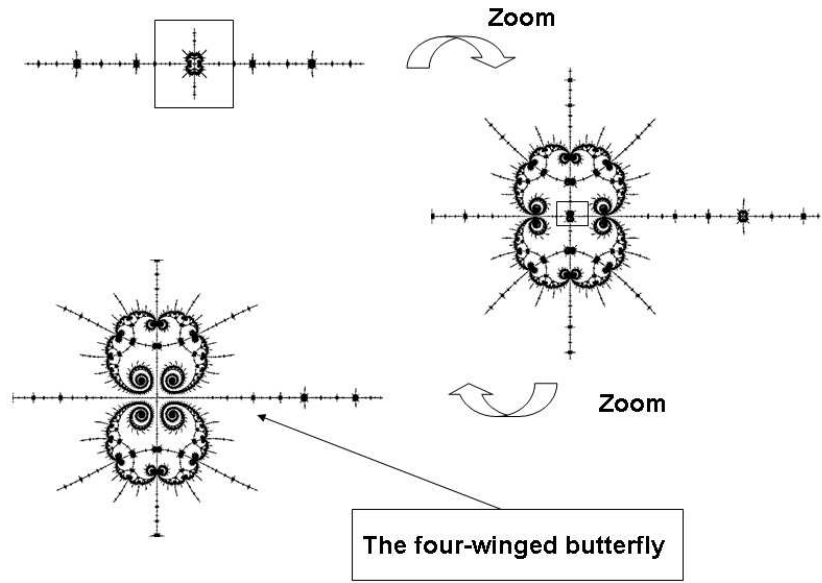}
\caption{The corresponding Julia set}
\end{center}
\end{figure}
\begin{figure}[!ht]
\begin{center}
\includegraphics[width=10cm, angle=270]{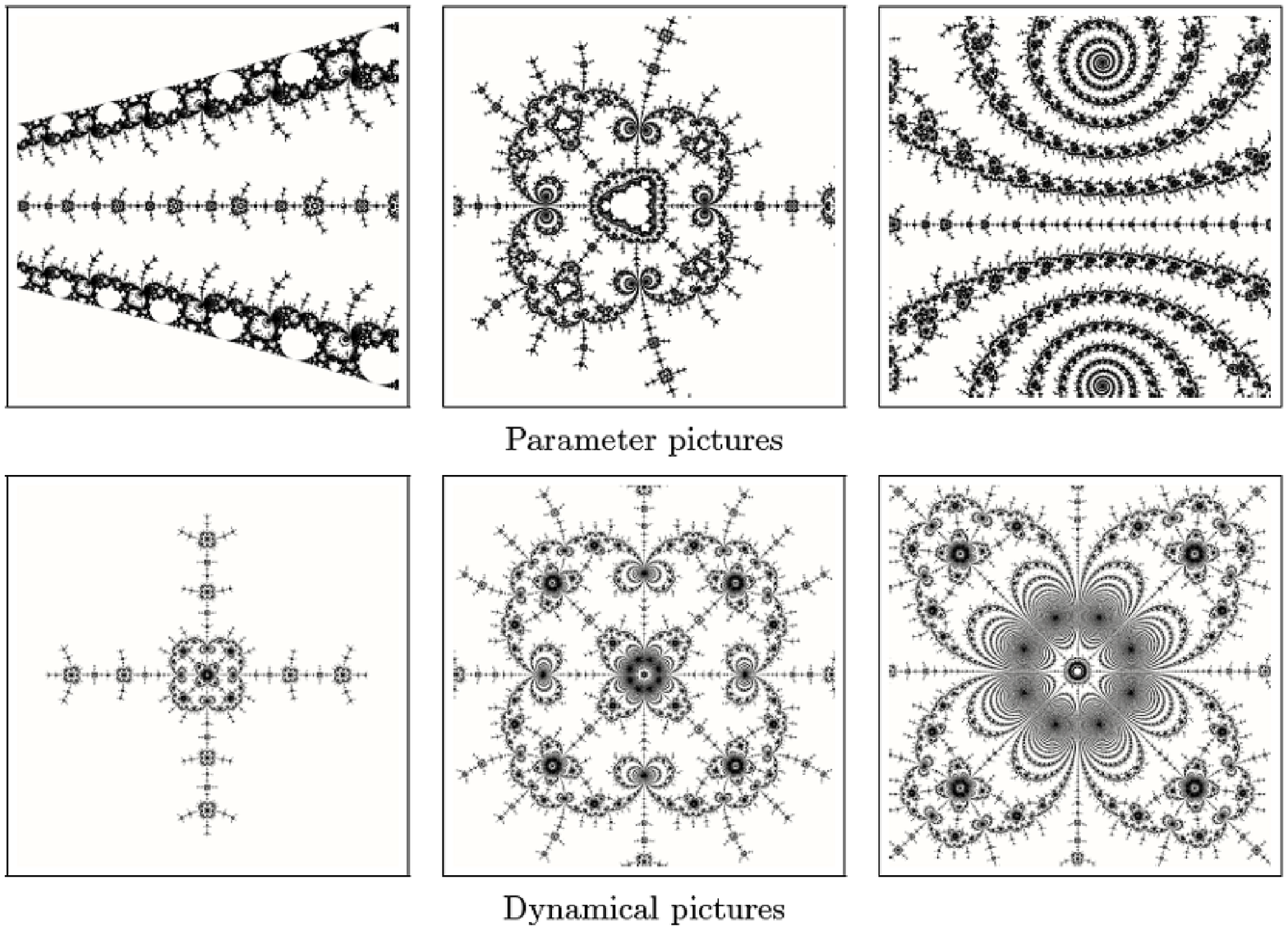}
\caption{The case d=4}
\end{center}
\end{figure}

Since $c_0\in (a, b)$, the set $J_c$ is connected for
$c$ in a small neighborhood $U$ of $c_0$. It is known,
that $f_c$ has an attracting periodic orbit of period $k$ for $c\in U$ on 
the left side of $c_0$.
Since $k>1$, the corresponding filled-in Julia set $K_{c_0}$ 
is such that 
its interior has a component $\Delta$ containing $0$, 
and infinitely many preimages of
$\Delta$. The boundary
$\partial \Delta$ contains a parabolic point $\alpha$ of period $k$.
In particular, there is a sequence of preimages of $\Delta$,
which intersect the real line and accumulate at $\alpha$.

For definicity, one can assume that $\alpha>0$.
Then $F=f_{c_0}^k$ 
has the following local form:
\begin{equation}\label{loc}
F(z)=z+a(z-\alpha)^2+b(z-\alpha)^3+...,
\end{equation}
where $a>0$.
This implies that a parabolic implosion phenomenon occurs as 
$\epsilon\to 0$, $\epsilon>0$, for the maps $f_{c_0+\epsilon}$.
\\ At this point we digress somewhat and describe briefly the theory of 
parabolic implosion.\\
Let $\delta>0$ be very small and $D_{\pm }$ be the disks centered at 
$\alpha\pm \delta$ with radius $\delta$. The map $F$ sends $D_-$ 
into itself while $D_+\subset F(D_+)$. This defines, after identification 
of $z$ with $F(z)$ at the boundary, two cylinders $U_-=D_-\backslash F(D_-)$ 
and $U_+=F(D_+)\backslash D_+$. 
The fact that $U_\pm $ are actual cylinders is best seen in the approximate 
Fatou coordinate $I: z\mapsto -1/(a(z-\alpha))$ which sends $\alpha$ to 
$\infty$ and conjugates 
$F$ to a map which is asymptotically the translation by $1$ at $\infty$:
\begin{equation}\label{Finf}
F_\infty(w)=I\circ F\circ I^{-1}(w)=w+1+\frac{A}{w}+O(\frac{1}{|w|^2}),
\end{equation}
where $A=1-b/a^2$. The real number $A$ is a conformal invariant. 
In the case of real polynomial $F$
which has a parabolic fixed point $\alpha$ with multiplier $1$ and 
with a single critical point in 
its immediate basin $\Delta$,  
it is known~\cite{Shi}
that $A>0$.\\
\begin{figure}[!ht]
\begin{center}
\includegraphics[width=9cm,angle=270]{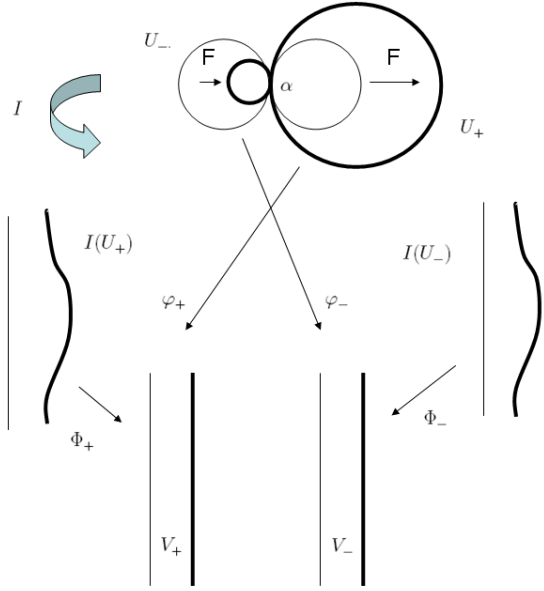}
\caption{Fatou coordinates}
\end{center}
\end{figure}

By Riemann mapping theorem these two cylinders may be uniformized 
by "straight" cylinders. In other words there exists $\varphi_\pm $ 
mapping the cylinders $U_\pm $ to vertical strips $V_\pm $ of width $1$ 
conjugating $F$ to the translation by $1$. For further use 
we notice that, by symmetry, $c_0$ being real, 
we may assume that $\varphi_\pm (\overline{z})=\overline{\varphi_\pm (z)}$.\\
We also notice that these maps are unique up to post-composition by a 
real translation. These two maps are called respectively repelling (+) 
and attracting (-) Fatou coordinates. We normalize them as follows.
For every $\kappa>0$ and $R>0$,
consider two sectors $\Sigma_-(\kappa, R)=\{w: Re(w)>R-\kappa |Im(w)|\}$, 
$\Sigma_+(\kappa, R)=\{w: Re(w)<-R+\kappa |Im(w)|\}$. 
Then for any $\kappa>0$ there is a big enough $R(\kappa)$,
such that, if we introduce two sectors $\Sigma_\pm(\kappa)=\Sigma_\pm(\kappa, R(\kappa))$, then
 $\varphi_\pm=\Phi_\pm\circ I$, where
\begin{equation}\label{phi}
\Phi_\pm(w)=w- A \log_\pm(w)+C_\pm+o(1)
\end{equation} 
as $w\to \infty$ within $\Sigma_\pm(\kappa)$ respectively.
We specify the constants $C_\pm$ and the log-branches in such a way, that $\Phi_\pm$ are real
for real $w$. Namely, we choose $C_-=0$ and $\log_-$ to be the standard log-branch in the slit plane
$\Ci\setminus \{x\le 0\}$. In turn, let $C_+=i A \pi$ and $\log_+$ to be a branch of $log$ in
$\Ci\setminus \{x\ge 0\}$ so that $\log_+(w)=\log|w|+i\pi$ for $w<0$.

Now, we extend $\varphi_\pm$ in the following way.
Since \begin{equation} \varphi_-(F(z))=T_{1} (\varphi_-(z)) \label{F-}\end{equation} where 
$T_\sigma$ denotes the translation by $\sigma$, and since every orbit converging to $\alpha$ 
passes through $U_-$ exactly once, $\varphi_-$ extends uniquely to $\Delta$ to an holomorphic 
function still satisfying ($\ref{F-}$). 


It is seen from~(\ref{F-}), that the map $\varphi_-: \Delta\to \Ci$ is a branched covering, with the
critical points at $0$ and all its preimages in $\Delta$ by $F^n$, $n>0$, and the critical values at
the real numbers $\varphi_-(0)-n$, $n\ge 0$. In particular, there exists a simply-connected domain 
$\Omega_-\subset \Delta\cap {\bf H}^+$, where ${\bf H}^+$ is the upper half-plane, such that
$\varphi_-: \Omega_-\to {\bf H}^+$ is a holomorphic homeomorphism. Moreover, the intersecrion of the
boundary of $\Omega_-$ with $\bf{R}$ is the interval $(0, \alpha)$.

Concerning the repelling Fatou coordinate it is best to consider $\psi_+= \varphi_-^{-1}:\, V_+\to U_+$. 
The functionnal relation is now \begin{equation}\psi_+(T_1(z))=F(\psi_+(z)) \label{F+}\end{equation} 
and we can extend $\psi_+$ to an entire function by putting, for $n\in\Z,\,\psi_+(T_n(z))=F^n(\psi_+(z))$ 
for $z\in V_+$. There exists a simply connected domain $\Omega_+\subset {\bf H}^+$, such that
$\psi_+: \Omega_+\to {\bf H}^+$ is a homeomorphism.\\
Let now $\sigma$ be a real number.
We define the Lavaurs map $g_\sigma$ on the component $\Delta$ of the
interior of the filled-in Julia set of $f_{c_0}$ 
by $$g_\sigma=\psi_+\circ T_\sigma\circ \varphi_-.$$
The "raison d'\^etre" of this definition is the following theorem due to Douady and Lavaurs ($\cite{AD}$):
\begin{theo}:There exists a sequence of positive
$\epsilon_n$
converging to $0$ and a sequence
of positive integers $N_n$ such that 
\begin{equation} g_\sigma(z)=\lim_{n\to \infty} f_{c_0+\epsilon_n}^{k N_n}(z) \label{Lavaurs}\end{equation}
uniformly on compact sets of $\Delta$. 
\end{theo}
Using~(\ref{phi}) with the constants $C_\pm$
and the log-branches specified as above, it is easy to get, 
that, for every $\kappa$,
if $w$ tends to $\infty$ in $\Sigma(\kappa):=\Sigma_-(\kappa)\cap \Sigma_+(\kappa)\cap {\bf H}^+$, then
$g_\infty(w):=I\circ g_\sigma\circ I^{-1}(w)=w+(\sigma-i A\pi)+O(\frac{1}{|w|})$,
where $A$ is real and positive.
Therefore, for every real $\sigma$ and every $\kappa>\kappa(\sigma)$, the inverse map
$g_\infty^{-1}$ leaves the sector $\Sigma(\kappa)$ invariant and $g_\infty^{-n}\to \infty$
as $n\to \infty$. Coming back to the $z$-plane, we conclude that the branch
$G=I^{-1}\circ g_\infty^{-1}\circ I$ of $g_\sigma^{-1}$ leaves the set $S(\kappa)=I^{-1}(\Sigma(\kappa))$
invariant, and 
$G^n(z)\to \alpha$
as $n\to \infty$, for $z\in S(\kappa)$ and every $\kappa>\kappa(\sigma)$. We have, for $w\in \Sigma(\kappa)$:
\begin{equation}\label{Ginf}
G_\infty(w):=I\circ G\circ I^{-1}(w)=w+(-\sigma+i A\pi)+O(\frac{1}{|w|}),
\end{equation}


Now, from the definition of $g_\sigma(z)$ and the global properties of the maps $\varphi_-$ and $\psi_+$, 
it follows the 
existence of a simply-connected domain $\Omega_0\subset \Omega_-$, which is mapped
by $g_\sigma$ homeomorphically onto ${\bf H}^+$ and such that $\alpha\in \bar \Omega_0$.
Moreover, from the above description, $S(\kappa)\subset \Omega_0$, for every $\kappa>\kappa(\sigma)$.
Therefore, the branch $G$ of $g^{-1}$ which is defined above, extends to a global univalent branch
$G: {\bf H}^+\to \Omega_0$ of $g_\sigma^{-1}$. 
Since $\Omega_0\subset {\bf H}^+$, the iterates $G^n: {\bf H}^+\to \Omega_0$, $n>0$, 
converge uniformly on compact sets in ${\bf H}^+$ to a unique fixed point in $\bar \Omega_0$,
which must be $\alpha$.

Let us consider the continuous map 
$\sigma\mapsto g_\sigma(0)=\psi_+(\varphi_-(0)+\sigma)$: 
if $\sigma$ runs in the interval 
$I=\{\varphi_+(x)-\varphi_-(0): x\in U_+\cap \R\}$, 
then $g_\sigma(0)$ runs over $U_+\cap \R$. It is thus clear, 
and this is the key point in the proof, 
that we can choose $\sigma$ in such a way that 
$g_\sigma(0)$ is a preimage of $\alpha$: there is $j\ge 1$, 
such that $f^j_{c_0}\circ g_\sigma(0)=\alpha$. 
Since $0$ is a critical point for $g_\sigma$,
taking the inverse image by $f^j_{c_0}\circ g_\sigma$ 
has the same effect as multiplying the number of petals by $d$
and we may state:
\begin{lem}\label{l}
There exists an infinite iterated function system
defined on a small compact neighborhood $B_0$ of zero and 
generated by some holomorphic branches of $f_{c_0}^{-1}$ and $g_{\sigma}^{-1}$
such that its limit set has Hausdorff dimension bigger than $2d/(d+1)$.
\end{lem}
\begin{com}\label{jl}
In fact, the limit set is a subset of a so-called Julia-Lavaurs set 
denoted by $J_{c_0, \sigma}$.
It is defined as follows.
The map $g_\sigma$ can be extended in a natural way from $\Delta$ to the interior
of the filled-in Fatou set of $f_{c_0}$: if $f_{c_0}^k(z)\in \Delta$, set
$g_\sigma(z)=g_\sigma\circ f_{c_0}^k(z)$. 
Then $J_{c_0, \sigma}$ is simply the closure of the set of 
points $z$ for which there exists $m\in\N$ such that $g_\sigma^m(z)$ is defined and belongs to $J(f_0)$.
\end{com}
This lemma together with the above discussion implies the theorem. Indeed,
by a general property of iterated function systems~\cite{mu}, 
there exists its finite subsystem with
the Hausdorff dimension of its limit set bigger than $2d/(d+1)$. On the other hand,
the finite iterated function system 
persists for $f_{c_0+\epsilon}$ by (\ref{Lavaurs}). To be more precise, if
$\{I_j: B_0\to X_j, 1\le j\le j_0\}$ is this finite iterated function system, then
each $I_j$ can be extended to a univelent map to a fixed neighborhood $Y$ of $B_0$ as $I_j: Y\to Y_j$.
Consider the inverse univalent map $I^{-1}_j: Y_j\to Y$. Since the convergence in (\ref{Lavaurs})
is uniform on compacts in $\Delta$, for every $\epsilon_n$ small enough
there is some integer $N_j>0$ and a compact set $X_{j, n}$, so that
$f_{c_0+\epsilon_n}^{N_j}: X_{j, n}\to B_0$ is univalent, too. Now it is clear, that, for
every $\epsilon_n$ small enough,
the non-escaping set $K_n$ of the dynamical system which consists of a finitely many maps
$f_{c_0+\epsilon_n}^{N_j}: X_{j, n}\to B_0$, $1\le j\le j_0$, has the Hausdorff dimension which is 
bigger than $2d/(d+1)$.     
On the other hand, $K_n$  
must lie in the Julia set
of $f_{c_0+\epsilon_n}$ because some iterate of the map
$f_{c_0+\epsilon_n}^{N_1 N_2 ... N_{j_0}}$
leaves the set $K_n$ invariant
and is expanding on it.
\paragraph{Proof of Lemma~\ref{l}.}
As the first step, let us fix a small enough closed ball $B_0$ around zero, so that it does not contain
points of the postcritical set of $f_{c_0}$.
There exists its preimage
$B'$ by $F^{-1}$ in $\Delta\cap {\bf H}^+$. 
Then we can apply to $B'$ the maps $G^n$, $n>0$. By the above, 
$B'_n=G^n(B')$ are pairwise disjoint, compactly contained in $\Delta$, and $B_n'\to \alpha$ as
$n\to \infty$. Now, for every $n\ge n_0$, so that $B_n'$ lies in a small enough neighborhood $U$ of $\alpha$,
we make ``clones'' of $B_n'$ in $U\cap \Delta$
applying to it $F^r$, $r\in\Z$, where $F^{r}$ for $r<0$
is a well-defined in $U\cap \Delta$
branch which fixes $\alpha$ . We obtain the sets $B'_{n,r}=F^r(B_n')$. 
On the second step, we consider the map $f_{c_0}^j\circ g_\sigma$ from a neighborhood $V$ of $0$ onto $U$,
This map is a remified cover with the only remification point at $0$ of order $d$.
Let $U^*=U\setminus \{x\ge \alpha\}$, and $V^*=V\cap \{z: Arg(z)\in (0, 2\pi/d)\}$. Denote 
by $h$ a branch of $(f_{c_0}^j\circ g_\sigma)^{-1}$ from $U^*$ onto $V^*$.
Let $B_{n,r}=h(B'_{n,r})$. We obtain a system of holomorphic maps
${\bf \Psi}=\{\psi_{n,r}: B_0\to B_{n,r}\}$, where 
$\psi_{n,r}=h\circ F^r\circ G^n\circ F^{-1}$, $r\in\Z, n\ge n_0$.
If the neighborhood $U$ is chosen small enough,
the maps $\psi_{n,r}$ extend to univalent maps in a fixed neighborhood $\tilde B$ of $B_0$
into itself. In particular, the compact sets $B_{n,r}$ are pairwise disjoint and compactly contained 
in $B_0$.
Now, it is quite standard to check that ${\bf \Psi}$ form a conformal infinite iterated function system
in the sense of~\cite{mu} (strictly speaking, in the hyperbolic metric of $\tilde B$, which is 
equivalent to the Eucledian one on $B_0$ though).
Let us calculate the parameter $\theta=\inf\{t: \psi(t)<\infty\}$ of ${\bf \Psi}$, where
$\psi(t)=\sum_{(n,r)} \max_{z\in B_0}|\psi_{n,r}'(z)|^t$. 
The map $\psi_{n,r}=h\circ I^{-1}\circ F_\infty^r\circ G_\infty^n\circ I\circ F^{-1}$.
Here $F^{-1}$ is a univalent map of a neighborhood $\tilde B$ of $B_0$ into $\Delta$.
Now, routine and well-know calculations based on~(\ref{Finf})-(\ref{Ginf}) show (see e.g.~\cite{Z}),
that, for all $r\in\Z, n\ge n_0$
and some $C$, which depends only on a compact set in ${\bf H}^+\cap \Delta$, from which $w$ is taken,
$C^{-1}|r+(-\sigma+i\pi A)n|\le |F_\infty^r\circ G_\infty^n(w)|\le C|r+(-\sigma+i\pi A)n|$,
and $C^{-1}\le |(F_\infty^r\circ G_\infty^n)'(w)|\le C$. On the other hand, the map $h$
is a composition of a univalent map with an inverse branch of $z^{1/d}$. This gives us:
$C_1^{-1} |r+(-\sigma+i\pi A)n|^{-1-1/d}\le |\psi_{n,r}'(z)|\le C_1 |r+(-\sigma+i\pi A)n|^{-1-1/d}$,
for some $C_1$ and every $z\in B_0$. It follows that the series for $\psi(t)$ converges if and only if
$t>\theta=2d/(d+1)$, and $\psi(\theta)=\infty$. Hence~\cite{mu}, the Hausdorff dimension
of the limit set of ${\bf \Psi}$ is strictly bigger than $2d/(d+1)$.



\small{Inst.\ of Math., Hebrew University, Jerusalem 91904, Israel,\\MAPMO, Universit\'e d'Orl\'eans, BP 6759 45067 Orl\'eans Cedex, France}\\

\bibliographystyle{plain}
\bibliography{biblio}

\end{document}